\documentclass[12pt]{article}
\usepackage{amsfonts} \usepackage{amssymb} \usepackage{amsthm}
\usepackage{a4}

\def\({\left(}
\def\){\right)}

\def\dim{{\rm dim\,}}

\def\tr{{\rm tr\, }}
\def\ad{{\rm ad\, }}

\def\eqref#1{(\ref{#1})}
\def\vol{\mathop\mathrm{vol}}
\def\area{\mathop\mathrm{area}}

\def\R{\mathbb{R}}
\def\G{\mathfrak{G}}
\def\N{\mathfrak{N}}

\def\proofend{\hspace*{1cm} \hfill $\rule{2mm}{2mm}$ \vspace{0.5cm}}
\def\semidirect{\rtimes}
\def\eqref#1{(\ref{#1})}
\def\ddt{\left.\frac{d}{dt}\right|_{t=0}}

\makeatletter
\def\section{\setcounter{equation}{0}\@startsection{section}{1}{\z@}%
                                   {-3.5ex \@plus -1ex \@minus -.2ex}%
                                   {2.3ex \@plus.2ex}%
                                   {\normalfont\large\bfseries}}
\def\subsection{\@startsection{subsection}{1}{\z@}%
                                   {-3.5ex \@plus -1ex \@minus -.2ex}%
                                   {2.3ex \@plus.2ex}%
                                   {\normalfont\bfseries}}
\makeatother

\newtheorem{thm}[equation]{Theorem}
\newtheorem{lm}[equation]{Lemma}

\newtheorem{coro}[equation]{Corollary}

\begin{document}

\title{The Cheeger constant of simply connected, solvable Lie groups}
\author{Norbert Peyerimhoff, Evangelia Samiou}
\date{}
\maketitle

\parbox{12cm}{
\small{\bf Abstract:}
We show that the Cheeger isoperimetric constant of a solvable simply
connected Lie group $G$ with Lie algebra $\G$ is
$h(G)=\max_{H\in\G,||H||=1} \tr(\ad (H))$.
}

\section{Introduction}
The Cheeger isoperimetric constant $h(M)$ of a complete noncompact
Riemannian manifold $M$ is defined by $$h(M)= \inf_{K\subset M}
\frac{\area (\partial K)}{\vol (K)},$$ where $K$ ranges over all
connected, open submanifolds of $M$ with compact closure and smooth
boundary.  For the bottom of the spectrum $\lambda(\Omega)$, there is
Cheeger's inequality \cite{c}, \cite {bpp}, \cite{cha}
$$\lambda(\Omega) \geq
\frac{h(\Omega)^2}{4}.$$
There is a converse of Cheeger's inequality due to Peter Buser,
\cite{bu},
\cite{bpp}: There exist constants $c_{1}$ and $c_{2}$, depending on a
lower bound on the Ricci curvature of $M$ such that
$$ \lambda\leq c_{1}h+c_{2}h^{2}.$$
Further relations between isometric and spectral
properties (in particular, estimates of the heat kernel) can be found
in the work of, e.g.,
%xxx
Varopoulos, Coulhon, Saloff-Coste, Grigor'yan \cite{gr}
and Pittet \cite{pi}. For connections between the Cheeger constant and
Kazhdan's property T see, e.g., \cite{Br} or \cite{Leu}.

%Cheeger constants are also useful in the theory of foliations. It was
%shown
%by Plante (\cite{pl}, \cite{su}) that if $M$ is a foliated Riemannian
%manifold
%with a leaf $\mathcal{F}\subset M$ satisfying $h(\mathcal{F})=0$ then
%the
%foliation admits an invariant transversal measure.

If $h(M)> 0$ one can easily see that $M$ has exponential volume
growth. The converse is not true. Hoke \cite{h} has shown that
$h(G)=0$ for a simply connected Lie group with left invariant metric
if and only if $G$ is unimodular and amenable. These Lie groups,
however, have exponential volume growth if they are not of type R, see
\cite{pa}. Examples of this class are horospheres in symmetric spaces
of noncompact type and higher rank orthogonal to the barycenter of a
Weyl chamber, see \cite{p}.

In this note we calculate the Cheeger constant for simply connected
solvable Lie groups.

\section{ The Cheeger constant}
We consider the following situation: Let $G=G_{0}\semidirect\R$ be a
semidirect product with Lie algebra
$\G = \G_{0} \oplus \R$ endowed with a left invariant metric such that
$\G_{0}$
is orthogonal to $\R$.
We denote the unit vector in $\R$ by $H_{0}$ and obtain
a diffeomorphism $\phi: G_{0} \times \R \to G$, $\phi(g_{0},t)=
g_{0} \exp(tH_{0})$.

\begin{lm}\label{mass}
The left invariant Haar-measure $\mu$ on $G$ is given by
\begin{equation}\label{measure}
d\mu(g_{0}, t) = \det(e^{\ad(-tH_{0})}) d\nu(g_{0})dt
= e^{-t\,\tr(\ad(H_{0}))} d\nu(g_{0})dt     \ ,
\end{equation}
where $d\nu$ is the left invariant Haar-measure on
$G_{0}=G_{0}\times\{0\}\subset G$.
\end{lm}
\noindent

Proof: Note that $d\nu dt$ is left $G_0$-invariant and right
$\R$-invariant.
Consequently, the left invariant Haar-measure on $G$ with respect to the

above
diffeomorphismus
is given by $d\mu(g_{0},t)=\delta(t)d\nu(g_{0})dt$. We will calculate
$\delta(t)$.
Let $a=\exp(sH_{0})$ and $f\in C^\infty_0(G) $. Then
\begin{eqnarray*}
\int_{G} f(ag)d\mu(g) &=& \int_{\R} \int_{G_{0}}
f(ag_{0}\exp(tH_{0}))\delta(t)d\nu(g_{0})dt  \\
                      &=& \int_{\R} \int_{G_{0}}
f(\psi(g_{0})a\exp(tH_{0}))\delta(t)d\nu(g_{0})dt,
\end{eqnarray*}
where $\psi:G_{0}\to G_{0}, \psi(g_{0})=ag_{0}a^{-1}$. From the
transformation formula
this becomes
$$ \int_{\R}\delta(t)\int_{G_{0}}f(g_{0} \exp((s+t)H_{0}))|\det
D\psi^{-1}(g_{0})|d\nu(g_0) dt\ .$$
Let $X_{1}, \cdots, X_{n}$ be an orthonormal basis of $\G_{0}$ with
respect to the left invariant metric.
Then the $Y_{j}=\ddt g_{0} \exp(tX_{j})$, $j=1,\ldots , n$ form an
orthonomal basis
of $T_{g_{0}}G_{0}$ and
\begin{eqnarray*}
D\psi^{-1}(g_{0})(Y_{j})&=&\ddt\psi^{-1}(g_{0}\exp(tX_{j}))
=\ddt a^{-1}g_0\exp(tX_j)a\\
&=& dL_{a^{-1}g_0 a}(e^{\ad(-sH_0)}X_j),
\end{eqnarray*}
so we conclude that $|\det(D\psi^{-1}(g_{0}))|=\det(e^{\ad(-sH_0)})$.
By left invariance of $\mu$ we get
\begin{eqnarray*} \int_{G} f(ag)d\mu(g) &=&
 \int_{\R}\delta(t)\det(e^{\ad(-sH_0)})
\int_{G_{0}}f(g_{0} \exp((s+t)H_{0}))d\nu(g_0) dt\ = \\
\int_{G} f(g)d\mu(g) &=&
\int_{\R}\delta(s+t) \int_{G_{0}} f(g_{0} \exp((s+t)H_{0}))d\nu(g_{0})dt
\ ,
\end{eqnarray*}
hence $\delta(s)= \det(e^{\ad(-sH_0)}).$
\proofend

The following theorem contains the main result of this section.
\begin{thm}\label{main}
For the Cheeger constant of a Lie group $G$ as above we have
\begin{equation} \label{mainest}
h(G) \geq |\tr(\ad(H_{0}))|
\end{equation}
Moreover, if $h(G_0) = 0$ then formula \eqref{mainest} holds with
equality.
\end{thm}

Proof: We identify $G$ with $G_0\times\R$ via the diffeomorphism
$\phi$ and denote by $\pi: G\to G_{0}$ the projection. Without loss of
generality we may assume that $\tr(\ad(H_0)) \ge 0$. For a connected,
compact subset $K\subset G$ with smooth boundary $\partial K$ and non
empty interior let $U=\pi(K)\setminus\{\mbox{critical values of }
\pi|_{\partial K}\}$.

On $U$ we define functions $\delta^\pm$ by
$\delta^+(u)=\max\{t\mid (u,t)\in K\}$
and $\delta^-(u)=\min\{t\mid (u,t)\in K\}$, $u\in U$.
The functions $\delta^\pm $ are smooth and by Sard's
theorem the set $U$ has full measure in $\pi(K)$. We denote by
$ \partial K^{\pm}=\{(u,\delta^{\pm}(u))| u\in U \} $
the graphs of $\delta^\pm$.
We estimate the volumes of $K$ and $\partial K$.
>From \eqref{measure} we get
\begin{eqnarray*}
\vol (K) & \leq & \int_{U}\int_{\delta^{-}(u)}^{\delta^{+}(u)}
e^{-t \,\tr(\ad( H_{0}))}dt d\nu(u)               \\
 & \leq & \frac{1}{\tr(\ad(H_{0}))}\int_{U} \(
e^{-\delta^{+}(u)\tr(\ad(H_{0}))}+
e^{-\delta^{-}(u)\tr(\ad(H_{0}))} \) \, d\nu(u) \ .
\end{eqnarray*}
Clearly $\area (\partial K)\geq \area (\partial K^+)+ \area (\partial
K^-)$.
We have
$$\area (\partial K^+)=
\int_U \sqrt{\det
(\langle
D\varphi(u)e_i,D\varphi(u)e_j\rangle_{(u,\delta^+(u))})_{i,j=1\ldots
n-1} }
\ d\nu(u),$$
where $\varphi(u)=(u,\delta^{+}(u))$ and $e_1,\ldots ,e_{n-1}$ is an
orthonormal basis of $T_uG_0$. We estimate the integrand
\begin{eqnarray}
\sqrt{\det(\langle
D\varphi(u)e_i,D\varphi(u)e_j\rangle_{(u,\delta^+(u))}) }
&=& \sqrt{\det(
\langle e_i,e_j\rangle_{(u,\delta^+(u))}+ vv^{\top}) } \nonumber \\
&\geq &\sqrt{\det\(\langle e_i,e_j\rangle_{(u,\delta^+(u))}\) }
\nonumber\\
&=&e^{-\delta^+(u)\tr(\ad(H_{0}))}, \label{norbi}
\end{eqnarray}
where $v^{\top}=(e_{1}(\delta^+),\cdots, e_{n-1}(\delta^+))$ and the
$e_i$'s
on the right hand side are considered as elements in
$T_{(u,\delta^{+}(u))}G$.
%xxx
The equality \eqref{norbi} follows from \eqref{measure}. This together
with the analogous estimate for $\partial K^-$ yields
\begin{equation} \label{form1}
\area (\partial K) \geq \tr(\ad(H_{0})) \vol (K) \ .
\end{equation}
This proves inequality \eqref{mainest}.

Finally, we prove equality of \eqref{mainest} in the case $h(G_0) =
0$. As before, we use the diffeomorphism $G \cong G_0\times \R$.
Let $K_0\subset G_0$ be arbitrary and consider the set
$K=K_0\times [a,b] \subset G$ with boundary
$\partial K= (K_0\times\{a,b\}) \cup (\partial K_0\times [a,b])$.
Direct calculations yield
\begin{equation}\label{volK}
\vol (K) = \vol (K_0) \frac{e^{- a \tr \ad H_0 }-e^{- b \tr \ad H_0 }}
{\tr\ad H_0} \ ,
\end{equation}
as well as
\[ \area(K_0\times\{a,b\}) = \vol(K_0) \left( e^{- a \tr \ad H_0} +
e^{- b \tr \ad H_0 } \right). \]
Given $\epsilon > 0$ we can arrange for
\begin{equation}\label{dk}
\frac{\area(K_0\times\{ a,b\})}{\vol (K)}\leq \tr \ad H_0 + \epsilon
\end{equation}
by choosing $b$ sufficiently large.
Note that \eqref{volK} and the estimate \eqref{dk} also make
sense if $\tr\ad H_0 =0 $.

Let $n = \dim \G_0$. Choosing an orthonormal basis $e_1,\ldots,
e_{n-1}$ of $T_u \partial K_0 \subset \G_0$ and abbreviating $e^{-t
(\ad H_0 )}=A_t$ we finally compute
$$\area(\partial K_0\times [a,b])=\int_a^b
\vol\nolimits_{n-1}(\partial K_0\times \{ t \})  \ dt$$
where the $(n-1)$-dimensional volume is given by
\[ \vol\nolimits_{n-1}(\partial K_0\times {t}) =
\int_{\partial K_0}\sqrt{\det\left(
\langle A_t e_p ,A_t e_q\rangle_{\G_0}\right)_{p,q=1\ldots n-1} }
\ d\vol\nolimits_{\partial K_0}(u).
\]
The integrand can be estimated from above by
$$M:=\max\left\{\left. \sqrt{\det\left(\left. P_U {A_t}^*A_{t}
\right|_{U}\right)} \ \right| \  a\leq t\leq b, \   U\subset \G_0, \
\dim U=n-1 \right\}, $$
where $P_U: \G_0 \to U$ denotes orthogonal projection. Note that $M$ is
independent of the specific choice of $K_0$. Therefore choosing $K_0$
such
that $\vol_{n-1}\partial (K_0)/\vol (K_0) $ is sufficiently small we may

assume that
\begin{equation} \label{form2}
\frac{\area(\partial K_0\times [a,b])}{\vol (K)}\leq \frac{M (\tr \ad
H_0)
(b-a)}{ e^{- a\tr \ad H_0 }-e^{- b\tr \ad H_0}}
\frac{\vol_{n-1}(\partial K_0)}{\vol (K_0)} \leq \epsilon.
\end{equation}
Putting \eqref{dk} and \eqref{form2} together yields $h(G) = \tr \ad
H_0$, since $\epsilon > 0$ was chosen arbitrarily small.
\proofend

Now we consider the particular case of a simply connected solvable Lie
group
$G$ with Lie algebra $\G$ and $\N=[\G, \G]$. In this case we obtain
\begin{coro}\label{cheegformel}
The Cheeger constant of a simply connected solvable Lie group $G$ with
Lie algebra $\G$ is
\begin{equation} \label{mainformula}
 h(G) = \max_{H\in \G, ||H||=1} \tr(\ad(H)) \ .
\end{equation}
\end{coro}
Proof: Let $\G_0$ denote the kernel of the $1$-form
$\alpha\colon\G\to\R$, $X \mapsto \tr\ad X $. Clearly $\G_0$ contains
$\N=[\G, \G]$, since $\N$ is nilpotent. Hence $\G_0$ is an ideal and the

corresponding unimodular Lie group satisfies $h(G_0) = 0$ by Theorem
\ref{main}, which proves the corollary in the case $\G_0 =
\G$. Otherwise $\G_0$ has codimension $1$. Let $H_0$ be the maximum of
$\alpha$ on the unit sphere of $\G$, i.e.
$$\max_{H\in \G, ||H||=1} \tr(\ad(H)) = \tr(\ad(H_0)) \ . $$
Then $H_0 \perp \G_0$ and, again, Theorem \ref{main} implies the
statement
of the Corollary.
\proofend

Remarks: Note that simply connected, solvable Lie groups with left
invariant metric may have curvature of both signs (e.g. horospheres in
symmetric spaces). Corollary \ref{cheegformel} in the particular case
of simply connected strictly negatively curved homogeneous spaces
(NCHS) was proved by Connell \cite{co}. He also showed for NCHS that the

Cheeger constant coincides with the exponential volume growth rate.

Formula \eqref{mainformula} also gives a lower bound for the topological
entropy
of compact Riemannian manifolds whose universal covering is a solvable
Lie group with left invariant metric. This follows from
the lower estimate of the topological entropy by the exponential
volume growth rate proved in \cite{ma}. In the case of a locally
symmetric space $M$ the topological entropy was calculated in
\cite{sp} and agrees with both the exponential volume growth rate and
the Cheeger constant of the universal covering $\tilde M$. It is given
by $\Vert \rho \Vert$, where $\rho$ is the sum of the positive roots
with multiplicities. More refined volume growth calculations in
symmetric spaces were carried out in \cite{kn}.

%It would be interesting to investigate the exponential volume growth
%rate of arbitrary solvable groups.

\bigskip

\noindent {\em Acknowledgements } We are grateful to P. Pansu for
bringing to our attention the paper of H. Hoke. We thank the University
of Cyprus for financial support.

\noindent
N. Peyerimhoff\\
Mathematische Fakult\"at\\
Ruhr-Universit\"at Bochum\\
Universit\"atsstr. 150\\
44780 Bochum\\
Germany\\
e-mail: peyerim@math.ruhr-uni-bochum.de

\bigskip

\noindent
E. Samiou\\
Department of Mathematics and Statistics\\
University of Cyprus\\
P.O. Box 20537\\
1678 Nicosia\\
Cyprus\\
e-mail: samiou@ucy.ac.cy

\end{document}